
\mag=1000 \hsize=6.8 true in \vsize=8.7 true in
 \baselineskip=15pt
\vglue 1.5 cm

       


 \def\qed{$\rlap{$\sqcap$}\sqcup$}


\font\tengothic=eufm10
\font\sevengothic=eufm7
\newfam\gothicfam
     \textfont\gothicfam=\tengothic
     \scriptfont\gothicfam=\sevengothic


  \font\tenmsb=msbm10              \font\sevenmsb=msbm7
\newfam\msbfam
     \textfont\msbfam=\tenmsb
     \scriptfont\msbfam=\sevenmsb
\def\Bbb#1{{\fam\msbfam #1}}



\def\move-in{\parshape=1.75true in 5true in}


\hyphenation {Castel-nuovo}

\def\twomed{\medskip\medskip}


\def\PP#1{{\Bbb P}^{#1}}
 
\def\P {\Bbb P}
\def\Y {\Bbb Y}
\def\X {\Bbb X}



\def\ref#1{[{\bf #1}]}

\def\Proof{\noindent {\it Proof:}}
\def\Prop#1{\noindent {\bf Proposition #1:}}
\def\Thm#1{\noindent {\bf Theorem #1:}}
\def\Lemma#1{\noindent {\bf Lemma #1:}}
\def\Coroll#1{\noindent {\bf Corollary #1:}}




\def\sec{\sigma_s(\X)}
\def\Xnd{\X_{{\underline n}, {\underline d}}}
\vskip 1cm \centerline{\bf On the ideals of Secant Varieties to
certain rational varieties.}
 
\bigskip\bigskip

\centerline{\it M.V.Catalisano, A.V.Geramita, A.Gimigliano.}
\bigskip
\bigskip

\move-in\noindent{\it Abstract:}\ If $\X \subset \P^n$ is a reduced and irreducible projective variety, it is interesting to find the equations describing the (higher) secant varieties of $\X$.  In this paper we find those equations in the following cases:

 \move-in$\bullet$ $\X = \P^{n_1}\times\ldots\times\P^{n_t}\times\P^n$ is the Segre embedding of the product and $n$ is 

\move-in``large" with respect to the $n_i$ (Theorem 2.4);

 \move-in $\bullet$ $\X$ is a Segre-Veronese embedding of some products with 2 or three factors;

\move-in $\bullet$ $\X$ is a Del Pezzo surface.

\bigskip\bigskip

\noindent {\bf 0. Introduction.}
\medskip

The study of the higher secant varieties of the Segre varieties has a long and interesting history (see e.g.\ref {ChCi}, \ref {ChCo}, \ref {K}, \ref {Pa}, \ref{Te}, \ref {Z} ). In addition to its intrinsic beauty
and its role in understanding properties of the projections of algebraic varieties, this study has been influenced by questions from representation theory, coding theory and algebraic complexity theory (see our paper \ref{CGG2} for some recent results as well as a summary of known results, and also \ref{BCS}).  Most surprising to us, however, are the connections with the recent work in algebraic statistics (e.g. see \ref {GHKM} and \ref {GSS}).

Although the major question classically asked about such secant varieties concerned their dimensions, and this is still - by and large - an open and challenging problem, the authors of the paper \ref {GSS} raised some interesting questions about the generators of the defining ideals of such varieties.

Unfortunately, questions about the commutative and homological algebra of the defining ideals of the higher secant varieties of any variety have received only limited attention.  Thus, apart from some
notable exceptions, there is very little information available about such questions.  One family of varieties for which we have rather complete information about the commutative algebra of their higher secant varieties is the family of rational normal curves (i.e. the {\it Veronese embeddings} of $\P^1$).  In this case the ideals in question are generated by the maximal minors of Hankel matrices and one knows not only these generators but also the entire minimal free resolution of these ideals.   Similarly, the
defining ideals for the higher secant varieties of the {\it quadratic Veronese embeddings} of $\P^n$ are defined by the (appropriately sized) minors of the generic symmetric matrix of size $(n+1)\times (n+1)$.  It follows, thanks to the work of [{\bf JPW}], that we thus know not only the generators of these ideals
but also their minimal free resolutions.

In this paper, however, our main interest is in Segre varieties.  In this case it is also well known that if the Segre variety is the embedding of
$$
\P ^{n_1} \times \P^{n_2} \longrightarrow \Y \subset \P^N, \ N =
(n_1+1)(n_2+1) -1
$$
and we let $\sigma_s(\Y)$ be the $(s-1)$-secant variety of $\Y$ (i.e. the closure of the union of all the $s$-secant $\P^{s-1}$'s to $\Y$) then $I_{\sigma_s(\Y)}$ is the ideal generated by the $(s+1)\times (s+1)$ minors of the $(n_1+1)\times (n_2+1)$ tensor (i.e. matrix) whose entries are the homogenous coordinates of
$\P^N$, i.e. the ideal of the $(s+1)\times(s+1)$ minors of the generic $(n_1+1)\times (n_2+1)$ matrix.  In this case the ideal is rather well understood (see e.g. {\ref L } and also the extensive bibliography given in the book of Weyman {\ref W}).

We will only refer to a small part of this vast subject and recall, e.g. that the ideal $I_{\sigma_s(\Y)}$ is a perfect ideal of height
$$
(n_1+1 - (s+1) -1)\times (n_2+1 - (s+1) - 1) = (n_1-s -1) \times (n_2-s-1)
$$
in the polynomial ring with $N +1$ variables, with a very well known resolution.

It follows from this description that all the secant varieties of the Segre embeddings of a product of {\it two} projective spaces are arithmetically Cohen-Macaulay varieties.  Moreover, from the resolution one can also deduce the degree, as well as other significant geometric invariants, of these varieties.

A determinantal formula for the degree was first given by Giambelli. There is, however, a reformulation of this result which we will use (see e.g.{\ref{H} pg. 244, or {\ref {BC, Thm. 6.5}}, where this lovely reformulation of the Giambelli Formula is attributed to J. Herzog and N.V. Trung):
$$
\deg (\sigma_s(\Y)) = \prod _{i=0}^{n_1 - s} {{{n_2+1+i\choose s} }
\over { {s+i\choose s} } }.
$$
It is worth mentioning that {\ref{BC, Thm. 6.9}} also have a very nice formula for the Hilbert Series of the coordinate ring of the various secant varieties to $\Y$, but we will not have occasion to use that formula here.

Let us now pass to the case of the Segre embeddings of {\bf more} than two factors.  More specifically, let $\X \subset \PP N$ denote the Segre embedding of
$$
\PP {n_1}\times \PP {n_2}\times ... \times \PP {n_t} \rightarrow \X \subset \PP N, \ N = \Pi^{t}_{i=1}(n_i+1)-1 \ , \ t \geq 3,
$$
where we usually assume that $n_1 \leq \cdots \leq n_t$.

If we let $T$ be the $(n_1+1)\times ...\times(n_t+1)$ tensor whose entries are the homogeneous coordinates in $\PP N$, then it is well known that the ideal of $\X$ is given by the $2\times 2$ minors of
$T$.  It is natural to ask if there is some way to use the tensor $T$ to get information about the higher secant varieties of $\X$.

If we partition $\{ 1, \ldots , t \}$ into two subsets (say $\{1, \ldots , \ell \}$ and $\{ \ell +1, \ldots , t\}$, to keep the notation simple) then we can form the composition
$$
(\P^{n_1}\times \cdots \times \P^{n_\ell}) \times (\P^{n_{\ell +1}} \times \cdots \times \P^t) \rightarrow \P^a \times \P^b
$$
where $a = \Pi _{i=1}^\ell (n_i +1) -1$, $b = \Pi_{i = \ell +1}^t (n_i + 1) - 1$, followed by
$$
\phi: \P^a \times \P^b \rightarrow \P^N, \ \ N \hbox{ as above}.
$$
Clearly $\phi(\P^a \times \P^b) \supseteq \X$ and hence
$$
\sigma_s(\X) \subseteq \sigma_s(\phi(\P^a \times \P^b)).
$$
Thus, the $(s+1)\times (s+1)$ minors of the matrix associated to the embedding $\phi$ will  all vanish on $\sigma_s(X)$.  That matrix, written in terms of the coordinates of the various $\P^{n_i}$ is called a {\it flattening} of the tensor $T$.

We can perform a {\it flattening} of $T$ for every partition of $\{ 1, \ldots , t \}$ into two subsets.  The $(s+1)\times (s+1)$ minors of all of these flattenings will give us equations which vanish on $\sigma_s(\X)$.  In \ref {GSS} it was conjectured that, at least for $s=2$, these equations are precisely the generators for the ideal $I_{\sigma_2(\X)}$ of $\sigma_2 (X)$.   The conjecture was recently proved in \ref {LM} for the special case of $t=3$ (and set theoretically for all $t$'s).  More recently, Allman and Rhodes (\ref{AR}) proved the conjecture for up to five factors, while Landsberg and Weyman (\ref{LW}) have found the generators for the defining ideals of secant varieties for the Segre varieties in the following cases: all secant varieties for $\P^1 \times \P^m \times \P^n$ for all $m, n$; the secant line varieties of the Segre varieties with four factors; the secant plane varieties for any
Segre variety with three factors. The proofs use representation theoretic methods.  To our knowledge, these are the only known results describing the ideals of higher secant variety for infinite families of Segre embeddings with more than 2 factors.

Note that for $s > 2$ one cannot expect, in general, that the ideals $I_{\sec}$ are generated by the $(s+1)\times (s+1)$ minors of flattenings of $T$.  Indeed, in many cases there are no such minors,
as the following example illustrates.

\medskip\noindent{\bf Example 0.1}\ Let $\X = \P^1 \times \P^1\times \P^1\times\P^1\times \P^1$ (5-times).
The Segre embedding gives us a $2\times 2\times 2\times 2\times 2$ box and the various flattenings will give us \medskip

\hskip 2cm $i)$\ 10 ($4\times 8$) matrices; and

\hskip 2cm $ii)$\ 5 ($2\times 16$) matrices.

\medskip
The largest minors we can look at for these flattenings are the $4\times 4$ minors of the first set of matrices and those will give us (some) equations for $\sigma_3(\X)$ and for no higher secant variety of $\X$.  But, $\X$ (of dimension 5) lives in $\P^{31}$ and one sees, by a simple dimension count, that $\sigma_4(\X)$ and $\sigma_5(\X)$ definitely lie on some hypersurfaces of $\P^{31}$.

\twomed Nevertheless, in the second section of this paper we will show that infinitely often the ideal of $\sec$ can be described by the $(s+1)\times (s+1)$ minors of {\bf ONE} flattening of $T$ (see Theorem 2.4 ). It follows immediately that these $\sigma_s(\X)$ are arithmetically Cohen-Macaulay schemes with a
well known minimal free resolution for their defining ideals.  As a consequence we obtain a method for finding the degrees of these secant varieties as well as other numerical invariants that can be calculated from the minimal free resolution (e.g. the Hilbert polynomial).

In the third section we study some Segre-Veronese varieties.  These are (special) linear sections of Segre varieties.   

In the final section of the paper we consider Del Pezzo varieties.  We give a complete description of the ideals of all of their secant varieties.

After this paper was written T. Abo, G. Ottaviani and C. Peterson obtained results about the {\it dimensions} of the secant varieties of products of projective spaces which are the same as our results on the dimensions in Theorem 2.4 (see \ref{AOP}).  We are very grateful to them for bringing this work to our attention. 

We would also like to thank E. Carlini for his help in making many of the computer calculations on which our conjectures and results are based.

\bigskip\noindent {\bf 1. Preliminaries.}

\medskip We will always work over an algebraically closed field $K$ of characteristic 0.

\medskip We recall the notion of {\it higher secant variety}.

\medskip\noindent {\bf Definition 1.1:} Let $\X\subseteq \PP N$ be a closed irreducible projective variety of dimension $n$. The $s^{th}$ {\it higher secant variety} of $\X$, denoted $\sec$, is the subvariety of $\P^N$ which is the closure of the union of all linear spaces spanned by $s$ linearly independent points of $\X$.

\medskip For $\X$ as above, a simple parameter count gives the following inequality involving the dimension of $\sec$ :
$$
\dim \sec \leq {\rm min} \{N, sn+s-1\}.  \eqno (1)
$$
Naturally, one ``expects" the inequality should, in general, be an equality.

When $\sec$ does not have the ``expected" dimension, $\X$ is said to be $(s-1)$-{\it defective}, and the positive integer
$$
\delta _{s-1}(\X) = {\rm min} \{N, sn+s-1\}-\dim \sec
$$
is called the  $(s-1)${\it-defect} of $\X$.

\medskip We will have occasion to consider a generalization of the higher secant varieties of a variety.  These are the {\it Grassmann secant varieties}, whose definition we now recall.

\medskip\noindent{\bf Definition 1.2:}\ Let $\X\subseteq \P ^N$ be a reduced and irreducible projective variety of dimension $n$, $s$ any integer $\leq N$.

For $k$ any integer, $0\leq k\leq s-1$, the {\it (k,s-1)-Grassmann secant variety of $\X$} (denoted $Sec_{k,s-1}(\X)$) is the Zariski closure, in the Grassmannian of $k$-dimensional linear subspaces of $\P ^N$ (which we will denote ${\Bbb G}(k,N)$) of the set
$$
\{l\in {\Bbb G}(k,N)\ \vert \ l\ \hbox{ is a subspace of the span of $s$ independent points of $\X$} \}.
$$
In case $k=0$ we get $Sec_{0,s-1}(\X) = \sec $.

\medskip As a generalization of the analogous result for the higher secant varieties, one always has the inequality
$$
\dim Sec_{k,s-1}(\X) \leq {\rm min}\{ sn+(k+1)(s-k-1), (k+1)(N-k)\},
$$
with equality being what is generally ``expected".

When  $Sec_{k,s-1}(\X)$ does not have the expected dimension then we say that $\X$ is $ (k,s-1)${\it-defective} and in this case we define the $(k,s-1)$-defect of $\X$ as the number:
$$
\delta _{k,s-1}(\X)= {\rm min}\{ sn+(k+1)(s-k-1), (k+1)(N-k)\}-\dim Sec_{k,s-1}(\X).
$$
(For general information on these defectivities see {\ref {ChCo}} and  \ref {DF}.)

In his paper \ref {Te2}, Terracini gives a link between these two kinds of defectivity for a variety $\X$ as above (see \ref {DF} for a modern proof):

\medskip \Prop {1.3} {\it (Terracini)} {\it Let $\X\subset \PP N$ be an irreducible non-degenerate projective variety of dimension n.  Let $\psi : \X \times \PP k \rightarrow \PP {(k+1)(N+1)-1}$ be the (usual) Segre embedding.

Then $\X$ is $(k,s-1)$-defective with defect $\delta_{k,s-1}(\X)=\delta$ if and only if $\psi (\X \times \PP k)$ is $(s-1)$-defective with $(s-1)$-defect $\delta_{s-1}(\X \times \PP k)=\delta$. }

\medskip\medskip Finally we wish to give a simple, but useful, lemma which we have been unable to find in the literature.

\medskip\noindent{\bf Lemma 1.4:} \ {\it Let $\X \subset \Y \subset \P^N$ be reduced irreducible projective varieties.  Suppose that for some integer $s$ we have:
$$
\sigma_s(\X) = \sigma_s(\Y)
$$
Then}
$$
\sigma_{s+1}(\X) = \sigma_{s+1}(\Y).
$$

\medskip\noindent{\it Proof:}\  One inclusion is clear, so suppose $P \in \sigma_{s+1}(\Y)$.  Then we can find $s+1$ linearly independent points of $\Y$, call them $Q_0, Q_1, \ldots , Q_s$, such that
$$
P = \alpha_0Q_0 + \alpha_1Q_1 + \cdots + \alpha_sQ_s
$$

Clearly $P^\prime = \alpha_1Q_1 + \cdots + \alpha_sQ_s \in \sigma_s(\Y) = \sigma_s(\X)$, so we can write
$$
P^\prime = \beta_1R_1 + \ldots + \beta_sR_s
$$
where $R_1, \ldots , R_s$ are linearly independent points in $\X$.  Thus, we can rewrite $P$ as
$$
P = \alpha_0Q_0 + \beta_1R_1 + \ldots + \beta_sR_s.
$$

Now consider
$$
P^{\prime\prime} = \alpha_0Q_0 + \beta_1R_1 + \ldots +
\beta_{s-1}R_{s-1} .
$$
With the same reasoning as above, we can write
$$
P^{\prime\prime} = \gamma_0T_0 + \ldots + \gamma_{s-1}T_{s-1}
$$
where $T_0, \ldots , T_{s-1}$are linearly independent points of $\X$.

Putting this all together we get
$$
P = \gamma_0T_0 + \ldots + \gamma_{s-1}T_{s-1} + \beta_sR_s
$$
and the points $T_0, \ldots , T_{s-1}, R_s$ are all points in $\X$.
That finishes the proof.

\bigskip\noindent {\bf 2. The main idea: the unbalanced case.}

\bigskip  As we mentioned earlier, we will be interested in finding Segre varieties
$\X$ for which some higher secant variety is described by the appropriate sized
minors of {\bf one} flattening of the tensor whose $2\times 2$ minors describe $\X$.  We will consider Segre embeddings of products
$$
\PP {n_1}\times \PP {n_2}\times ... \times \PP {n_t}\times \PP n,
$$
where $n_1\leq n_2\leq ... \leq n_t \leq n$ (often $n >> n_t$ hence the term ``{\it unbalanced}").

\medskip\medskip The following easy example (see \ref{P}, and for the case of the secant line variety see also \ref{LM}) will illustrate  the main idea in what follows.

\bigskip\noindent {\bf Example 2.1:} \ Consider the Segre varieties $\X$ given by embedding $\PP 1\times \PP 1\times \PP n$ into $\PP {4n+3}$, $n\geq 2$. The ideal of $\X$ is given by the $2\times 2$ minors of a $2\times2\times (n+1)$ tensor $T$ of indeterminates (the coordinates of $\PP {4n+3}$).
$$
T\ \ \ \ \ =\ \ \ \ \ \ \ \ \matrix{v_{000}& ---&v_{001}&---&v_{002}& ...&v_{00n}&&\cr
&&&&&&&&\cr |&\searrow&&\searrow&&&&\searrow&\cr &&&&&&&&\cr
v_{100}&&v_{010}&---&v_{011}&---&...&---&v_{01n}\cr &&&&&&&&\cr
&\searrow&|&&|&&&&|\cr
&&&&&&&&\cr&&v_{110}&---&v_{111}&---&...&---&v_{11n}}
$$
\medskip
Consider the $4\times (n+1)$ matrix $M$ obtained by flattening this tensor, i.e. by using the composition:
$$
(\P^1 \times \P^1) \times \P^n\ \rightarrow \  \P^3 \times \P^n \rightarrow \P^{4n+3},
$$
$$
M = \pmatrix {v_{000}& v_{001}& v_{002}& ---&v_{00n}\cr
  & & & &  \cr
v_{100}& v_{101}& v_{102}& ---&v_{10n}\cr
 & & & & \cr
v_{110}& v_{111}& v_{112}& ---&v_{11n}\cr &&&&\cr v_{010}&
v_{011}& v_{012}& ---&v_{01n}}.
$$
The ideal generated by the $2\times 2$ minors of $M$ is the ideal of the Segre variety $\Y$ given by  embedding $\PP 3\times \PP n$ into $\P^{4n+3}$.  Trivially $\Y$ contains $\X$.

Now consider the ideal generated by the $3\times 3$ minors of $M$.  This is well known to be the ideal of $\sigma_2(\Y)$, and of course those minors also vanish on $\sigma_2(\X)$. Since the matrix $M$ has
generic entries, we know that its $3\times 3$ minors generate a prime ideal of height $(4-3+1)(n+1-3+1)= 2n-2=$ codim $\sigma_2(\Y)$, which in fact is defective (its expected codimension is $2n-4$).  Thus the dimension of $\sigma_2(\Y)$ is $2n+5$.

But we know that $\dim \sigma_2(\X) = 2n+5$ (Segre varieties with three or more factors always have $\sigma_2$ of the expected dimension, see \ref{CGG2}), hence $\sigma_2(\X)=\sigma_2(\Y)$. By the Lemma above, this implies that (for $t\geq 3$) also $\sigma_t(\X)=\sigma_t(\Y)$.  In this example, the only other
relevant $t$ is $t = 3$ (for $t = 4$ we have $\sigma_4(\X) = \P^{4n+3}$) and thus $I_{\sigma_3(\X)}$ is generated by the $4\times 4$ minors of $M$ (only relevant when $n \geq 3$).  This ideal is an ideal of height $(n+1-4+1)=n-2$ and so $\sigma_3(\X)$ is defective, having dimension $3n+5$ instead of $3n+7$.  Since $\sigma_3(\X)$ is defined by the maximal minors of a $4 \times (n+1)$ matrix, we can also say that its degree is ${n+1\choose 3}$ and it is arithmetically Cohen-Macaulay.

\medskip  The point of this example is, we hope, clear: it sometimes happens that a $t^{th}$ secant variety for the Segre product of three or more projective spaces is the same as the $t^{th}$ secant variety of a Segre product with only two factors.  Inasmuch as we have abundant information about Segre products with two factors, this gives us a way to get information about Segre products with more than two factors.

\medskip\medskip  Our first task is to find more times when the behavior in Example 2.1 occurs.
This  is the content of the following Lemma.

\medskip\Lemma {2.2} {\it Let $V\subset \PP N$ be a variety such that $Sec_{s-1,s-1}(V) = {\Bbb G}(s-1,N)$.
Consider  the Segre embedding $\Y$ of $\PP N\times \PP n$ into $\PP M$, $M=Nn+N+n$. 

If $\X$ is the image of $V\times \PP n$ into $\PP M$, then $\sigma_s(\X)=\sigma_s(\Y)$.}

\medskip\noindent {\it Proof:} Let $\phi :\PP N\times \PP {n} \longrightarrow\PP M $ be the Segre embedding.   Consider a general secant $\P^{s-1}$ to $\Y$ (the image of $\phi$) and call it $H$.  Then,
$$
\PP {s-1}\cong H = <\phi(A_0,B_0),...,\phi(A_{s-1},B_{s-1})>,
$$
with $A_i \in \PP N$, $B_i \in \PP n$, generic points in their spaces. For all $\lambda_0,...,\lambda _{s-1}\in K$ we want to check that the point:
$$
 P_{\underline \lambda}= \lambda_0\phi(A_0,B_0)+...+\lambda_{s-1}\phi(A_{s-1},B_{s-1})
$$
is in $\sigma_s(\X)$.  We will be done if we find points $C_0,...,C_{s-1}$ in $V$ and  $D_0,...,D_{s-1}$ in $\PP n$  such that:
$$
P_{\underline \lambda}= \phi(C_0,D_0)+...+\phi(C_{s-1},D_{s-1}).
$$

Since $Sec_{s-1,s-1}(V) = {\Bbb G}(s-1,V)$ and the points $A_i$ are generic in $\P^N$ we can choose the $C_i$'s in $V$ such that
$$
<C_0,...,C_{s-1}>\ = \ <A_0,...,A_{s-1}>,
$$  
and so we can write
$$
A_i=\sum_{j=0}^{s-1}a_j^{(i)}C_j.
$$

Since $\phi$ is a bilinear map, we obtain:
$$ 
 P_{\underline \lambda}= \lambda_0\phi(A_0,B_0)+...+\lambda_{s-1}\phi(A_{s-1},B_{s-1})
$$
$$
=\lambda_0\phi(\sum_{j=0}^{s-1}a_j^{(0)}C_j,B_0)+...+\lambda_{s-1}\phi(\sum_{j=0}^{s-1}a_j^{(s-1)}C_j,B_{s-1})
$$
$$=\lambda_0 \left[ \sum_{j=0}^{s-1}a_j^{(0)}\phi (C_j,B_0)\right] +...+ \lambda_{s-1}\left[ \sum_{j=0}^{s-1}a_j^{(s-1)}\phi (C_j,B_{s-1})\right]
$$
$$
=\sum_{j=0}^{s-1}\phi (C_j,\lambda_0 a_j^{(0)}B_0)+...+ \sum_{j=0}^{s-1}\phi (C_j,\lambda_{s-1} a_j^{(s-1)} B_{s-1})
$$
$$
=\sum_{j=0}^{s-1}\phi (C_j,\lambda_0 a_j^{(0)}B_0+ ... + \lambda_{s-1} a_j^{(s-1)} B_{s-1})
$$
$$ 
 = \phi(C_0,D_0)+...+\phi(C_{s-1},D_{s-1})
$$
where
$D_j:=\lambda_0 a_j^{(0)}B_0+ ... + \lambda_{s-1} a_j^{(s-1)} B_{s-1}\in \PP n$ and the $C_i$'s are in $V$, and we are done. 

\hfill\qed

\bigskip \bigskip We now look for times when the hypothesis of Lemma 2.2 are satisfied.  To that end, consider Segre varieties $\X\subset \PP M$, $M= [\Pi_{i=1}^t(n_i+1)](n+1) - 1$, given by embedding  $\PP
{n_1}\times \PP {n_2}\times ... \times \PP {n_t}\times \PP n$ into $\P^M$, and also the Segre variety $\X'$ given by embedding $\PP {n_1}\times \PP {n_2}\times ... \times \PP {n_t}$ into $\PP N$,
$N=\Pi_{i=1}^t(n_i+1) - 1$.  We can consider $\X$ as obtained by composing the map, which is the Segre embedding on the first factor and the identity on the second factor, with the Segre embedding $\Y$
of $\PP N\times \PP n$ into $\PP M$.  I.e. we have:
$$
\matrix{\PP {n_1}\times \PP {n_2}\times ... \times \PP {n_t}\times
\PP n&&\longrightarrow&&\PP M \cr \downarrow &&&\nearrow & \cr
\X'\times\PP n&\subset&\PP N\times \PP n&&}
$$

\bigskip\noindent{\bf Lemma 2.3:}\ {\it The Segre variety $\X' = \PP {n_1}\times \PP {n_2}\times ... \times \PP {n_t}\subset \PP N$ has the property that $Sec_{s-1,s-1}(\X') = {\Bbb G}(s-1,N)$ if and only if:}
$$
N-\sum_{i=1}^tn_i +1\leq s .
$$
 
\medskip\noindent{\it Proof:}\  Since $\dim \X^\prime = \sum_{i=1}^t n_i$, if $s \leq N-\sum_{i=1}^tn_i = codim (\X^\prime)$, then a generic $\PP {s-1}\subset \PP N$ will not intersect $\X^\prime$, hence in this case
$Sec_{s-1,s-1}(\X') \neq {\Bbb G}(s-1,N)$.

Now let $N-\sum_{i=1}^tn_i +1= s$; since $\dim \X^\prime = \sum_{i=1}^t n_i$ and $\X^\prime$ is reduced and non-degenerate, a general linear subspace of $\P^N$ of dimension $N - \sum_{i=1}^t n_i$ will meet $\X^\prime$ in $\deg \X^\prime$ distinct points.  Again, since $\X^\prime$ is non-degenerate, $\deg \X^\prime \geq$ codim $\X^\prime = (N - \sum_{i=1}^t n_i) + 1$.  Thus, since $s-1 = N - \sum_{i=1}^t n_i$, a generic $\P^{s-1}$ of $\P^N$ meets $\X^\prime$ in at least $s$ points.  Hence, such a $\P^{s-1}$ is
definitely an $s$-secant linear space to $\X^\prime$.  It follows that for this $s$ we have
$$
Sec_{s-1,s-1}(\X^\prime) = {\Bbb G}(s-1, N) . \eqno{(*)}
$$

If now we choose $s$ so that $s-1 > N - \sum_{i=1}^t n_i$ then a generic $\P^{s-1}$ of $\P^N$ will meet $\X^\prime$ in a variety of dimension $> 0$ and hence will certainly be a secant $\P^{s-1}$ to $\X^\prime$.  Thus, $(*)$ is also true for such an $s$ and Lemma 2.3 has been verified.

\hfill\qed

\twomed With all the preliminary observations being established, we are now ready to prove the main result of this section.
 
\bigskip \Thm {2.4} {\it Let $\X$ be the Segre embedding
$$
 \P^{n_1}\times \ldots \times \P^{n_t} \times \P^n \rightarrow \X \subset \P^M, \ \ \ \ M =(n+1)(\Pi_{i=1}^t(n_i +1)) - 1 ,
$$ 
and let $\Y$ be the Segre embedding of $\P^N \times \P^n$ in $\P^M$, $N =  \Pi_{i=1}^t(n_i +1) - 1$.  Let $n\geq  N - \sum _{i=1}^tn_i + 1$. 

Then:

\medskip\move-in\noindent 1) \ for $2 \leq s \leq N-\sum_{i=1}^tn_i$,   \  
$\sigma_s (\X)\neq \sigma_s (\Y)$  and $\sigma_s (\X)$ has the expected dimension;

\medskip\move-in\noindent 2) for $ s = N-\sum_{i=1}^tn_i + 1$, \quad $\sigma_s
(\X) = \sigma_s(\Y)\neq \P^M$ and $\sigma_s(\X)$  has the expected dimension;

\medskip\move-in\noindent 3) for $N-\sum_{i=1}^tn_i +1 < s \leq \min\{n,N\}$,\quad
$\sigma_s (\X) = \sigma_s(\Y)\neq \PP M$ and $\sigma_s (\X)$ is
defective with $\delta_{s-1} (\X) = s^2-s(N-\sum_{i=1}^tn_i +1)$;

\medskip\move-in\noindent 4) for $s \geq \min\{n,N\}+1$, \quad  $\sigma_s (\X)=
\sigma_s (\Y)= \PP M$;

\medskip\move-in\noindent 5) in cases 2) and 3) above, the ideal of $\sigma_s (\X)=
\sigma_s (\Y)$ is generated by the $(s+1)\times (s+1)$ minors of an
$(n+1)\times (N+1)$ matrix of indeterminates.   

\medskip\move-in It follows that, in cases 2) and 3), $\sigma_s (\X)$ is a.C.M. and a minimal free resolution of its defining ideal is given by the Eagon-Northcott complex.
}

\medskip \Proof \ $2)$ First notice that from Lemmas 2.2 and 2.3, the equality $\sigma_s (\X)= \sigma_s (\Y)$ is immediate.  We have already mentioned that one knows the dimension of $\sigma_s(\Y)$ for any $s$ and a simple calculation reveals that the dimension we obtain for $\sigma_s(\X)$ is that which is expected.

As for $1)$, the hypothesis $n \geq N - \sum_{i=1}^tn_i +1$ guarantees that $\sigma_s(\Y) \neq \P^M$.  The result then follows immediately from $2)$ and our knowledge of the dimensions of $\sigma_s(\Y)$.

The equality of $\sigma_s(\X)$ and $\sigma_s(\Y)$ in $3)$ and $4)$ is again guaranteed by Lemmas 2.2 and  2.3.  Once again we use the fact that the dimensions of the $\sigma_s(\Y)$ are known and a simple calculation gives: the defectivity in the range described in $3)$; the equality in the range described in $4)$. 

$5)$ is, again, an immediate application of our characterization of the ideal of $\sigma_s(\Y)$. The closing statement of the theorem also follows from this characterization.  \hfill\qed

\twomed\noindent{\bf Remark 2.5:}\ If we continue with the notation of Theorem 2.4, and suppose that $n = N - \sum_{i = i}^tn_i$ and $s = n+1$ then $\sigma_s(\X) = \sigma_s(\Y) = \P^M$ (exactly as in part $4)$ of Theorem 2.4).  Moreover, one has the additional fact that
$$
\dim \sigma_s(\X) = s \dim(\X) + (s-1)
$$
and hence that $\dim \sigma_t(\X)$ is the expected dimension for {\it every} $t$.

\bigskip As a consequence of Theorem 2.4 we have the following:

\bigskip\noindent {\bf Corollary 2.6:}  {\it Let $\X\subset \PP M$ be the Segre embedding of $\PP 1\times
\PP m \times \PP n$, $m\leq n$ (hence $M=2nm+2n+2m+1$).  

\noindent 

\medskip\    i)\ If\  $n=m$, then $\sigma _s(\X)$ has the expected dimension for all $s$;

\medskip\    ii)\ if\  $n = m+1$, then $\sigma _s(\X)$ has the expected dimension for all $s$; 

\medskip\    iii)\ if\   $n > m+1$, and  

\medskip\hskip 20pt\noindent    $\bullet\  2 \leq s \leq m+1$, then $\sigma _s(\X)\neq \P^M$ has the expected dimension;  

\medskip\hskip 20pt\noindent    $\bullet\  m+2 \leq s \leq \min\{ 2m+1, n \}$, then $\sigma _s(\X)$  is defective with  $\delta_{s-1} (\X) = s^2-s(m+2)$;

\medskip\hskip 20pt\noindent     $\bullet\  s > \min \{2m+1,n\}$, then $\sigma _s(\X)= \PP M$.

}

\medskip\Proof  \ $i)$ is immediate from Remark 2.5.  $ii)$ and $iii)$ all follow from the various parts of Theorem 2.4.  

\hfill\qed

\bigskip
Notice that partial results in this case can be found in \ref{J}, see \ref{CGG2}, pg 282. It is interesting to compare our results with those found by
\ref{CS}, \ref{LM} and \ref{LW}.

\bigskip\noindent {\bf Example 2.7:} The family of Segre varieties
$\P^1 \times \P^m \times \P^n$, has been also considered in
\ref{LW}. These authors show that the ideal of $\sigma_s(\X)$ is
generated by the $(s+1)\times (s+1)$ minors of the flattenings of
the tensor $T$ giving the embedding of $\X$.
 \ref{LW} do not discuss the dimensions of the secant varieties to members of
this family and, consequently, do not mention their defectivities.
Note, however, that there are really only two flattenings to consider
for members of this family (the third one has no $3\times 3$
minors). In fact, for $m+1 \leq s \leq \min\{n, 2m+1 \}$,
$I_{\sigma_s(\X)}$ is the ideal of the $(s+1)\times (s+1)$ minors of
a single flattening of $T$. The proofs in \ref{LW} rely on a
subtle analysis using representation theory.

In any case, when we have that $\sigma_s(\X)$ is determinantal, then
it is given by the $(s+1)\times (s+1)$ minors of a single
flattening. We can then apply the Giambelli formula in order to get the degree of
$\sigma_{s}(\X)$. For example, if we consider the case $n = m+1$ and
let $\X_m$ be the Segre embedding of $\PP 1\times \PP m \times \PP {m+1}$, then
$\sigma_{m+1}(\X)$ has ideal generated by the $(m+2)\times (m+2)$
minors of a $2(m+1)\times (m+2)$ matrix and hence
$$
\deg (\sigma_{m+1}(\X_m)) = {2m+2 \choose m+1}.
$$
Cox and Sidman (in \ref{CS, Theorem 5.1}) give a formula for the
degree of $\sigma_2(\X_m)$.  For $s$ such that $3 \leq s \leq m$ we
are not aware of any method to calculate the degree of
$\sigma_s(\X_m)$.

\bigskip Now let us consider the case of four factors $ 
\P^{n_1}\times\P^{n_2}\times\P^{n_3}\times\P^{n_4} \rightarrow \X \subset \P^M, \ M
=  [ \prod_{i=1}^4(n_i + 1) ] -1$ (and, as always, $n_1 \leq n_2
\leq n_3 \leq n_4$).

In this case \ref{LW} prove that the ideal of $\sigma_2(\X)$ is
generated by the $3\times 3$ minors of all the flattenings of the
tensor describing $\X$.

If we consider the function
$$
N(n_1,n_2,n_3) = (n_1+1)(n_2+1)(n_3+1) - (n_1+n_2+n_3)
$$
then our results apply to all those $\X$ (as above) for which $n_4
\geq N(n_1,n_2,n_3)$.  In this case we have:

\move-in $1)$\ a complete description of the dimensions of
$\sigma_s(\X)$ for every $s$;

\medskip\move-in $2)$ \ if, in addition, $s \geq N(n_1,n_2,n_3)$
then the ideal of $\sigma_s(\X)$ is generated by the minors of one flattening of the tensor describing$\X$ and so we also know the finite free resolution of this ideal.  

\medskip These results apply, for example, to:
$$
\matrix{ \P^1\times \P^1 \times \P^1 \times \P^n & \hbox{ for } n
\geq 5; \cr
 \P^1\times \P^1 \times \P^2 \times \P^n & \hbox{ for } n \geq 8; \cr
 \P^1\times \P^2 \times \P^2 \times \P^n & \hbox{ for } n \geq 13; \cr
 \P^2\times \P^2 \times \P^2 \times \P^n & \hbox{ for } n \geq 21. }
$$

\bigskip It is also possible to apply Theorem 2.4  in order to obtain results on
Grassmann defectivity.

\bigskip
\Coroll {2.8} {\it Let $\X, n, N$ be as in Theorem 2.4, and let $\X
_i$ be the Segre embedding of $\PP {n_1}\times ...\times\widehat \PP
{n_i}\times...\times \PP {n_t}\times\PP {n}$, $i=1,...,t$. Then for
$N-\sum_{i=1}^t n_i+1\leq s\leq \min\{n,N\}$ and $s\geq n_i+1$, we
have that $Sec_{n_i,s-1}(\X_i)$ is defective, with
$\delta_{n_i,s-1}=s^2-s(N-\sum_{i=1}^t n_i+1)$, while
$Sec_{n_i,s-1}(\X_i)$ has the expected dimension for all other
values of $s$.}

\medskip\Proof \ The corollary is a direct consequence of Theorem
2.4 and Proposition 1.3.\hfill\qed

\bigskip

\Coroll {2.9} {\it Let $\X_n$ be the Segre embedding of $\PP 1\times
\PP n$. Then $Sec_{m,s-1}(\X_n)$ is defective if and only if $n>m+1$
and $m+2\leq s\leq \min\{2m+1,n\}$. Moreover, in this case,
$\delta_{m,s-1}(\X_n)= s^2-s(m+2)$.}

\medskip
\Proof \  This follows directly from Corollaries 2.6 and 2.8. \hfill\qed

\bigskip
In the following example we will consider how to use the \lq\lq unbalanced
case" idea also when we are not able to describe the ideal of the
secant variety completely.

\twomed\noindent{\bf Example 2.10:}\ In this example we would like to consider the following family of Segre varieties, this time with four factors: 
$$
  \P^1 \times \P^1 \times \P^n \times \P^n \rightarrow \X(n) \subset \P^N, \ \ \  N =
4n^2 + 8n + 3.
$$

Theorem 2.4 does not apply to members of this family.  Thus, we cannot say that any secant
variety for this family has equations derived from {\bf one}
flattening of the tensor describing $\X(n)$. Nevertheless, it is
possible to show, using flattenings, that $\sigma_{2n+1}(\X(n))$ is
defective.

A quick check shows that the ``expected dimension" of
$\sigma_{2n+1}(\X(n))$ is $4n^2 + 8n + 2$, i.e. we expect that
$\sigma_{2n+1}(\X(n))$ is a hypersurface of $\P^N$.  But, if we
group the factors of $\X(n)$ above as
$$
(\P^1\times \P^n) \times (\P^1\times \P^n)
$$
and then permute the $\P^1$'s, we obtain two distinct embeddings of
$\P^{2n+1}\times \P^{2n+1} \rightarrow \P^N$ and the determinants of
the resulting matrices of size $(2n+2)\times (2n+2)$ give us two linearly independent 
forms of degree $2n+2$ which vanish on $\X(n)$.  Consequently,
$$
\dim \sigma_{2n+1}(\X(n)) \leq N-2
$$
and hence it is defective.

In case $n = 1$, we showed in \ref{CGG3, Example 2.2} that $\dim
\sigma_3(\X(1)) = 13$.  This is precisely $N-2$ for this case.  One
can show that the ideal of $\sigma_{3}(\X(1))$ is generated by two
quartics (hence these two) even though for $n=1$ there is yet a
third flattening which gives a third quartic in the ideal.  But of
these three quartics, any two generate the ideal and the third is a
linear combination of the other two.

This raises several interesting questions for this family:

\move-in\noindent $1)$\ Is $\sigma_{2n+1}(\X(n))$ the only defective
secant variety for $\X(n)$?  From \ref{CGG2, Prop. 3.7} we know that
$\sigma_s(\X(n))$ is not defective for $s \leq n+1$.

\medskip\move-in\noindent $2)$\ Is $\sigma_{2n+1}(\X(n))$ always the
complete intersection of the two forms of degree $2n+2$ that we found above?

\twomed\noindent{\bf Remark 2.11:}\ Since this preprint was distributed, \ref{AOP} resolved the first question.  They showed that the codimension of $\sigma_{2n+1}(\X(n))$ is exactly two.  They also showed that it is the only defective secant variety in this family using their induction procedure.  In fact, this last follows immediately from the knowledge that the codimension is exactly two.

Our reasoning, which differs from that in \ref{AOP}, goes as follows: given the codimension, one knows that the defectivity of the varieties $\sigma_{2n+1}(\X(n))$ is exactly 1 and hence that the secant varieties $\sigma_t(\X(n)), \ t \leq 2n$ (which must have smaller defectivity) cannot be defective at all.  It is easy to check, using the fact that the codimension of $\sigma_{2n+1}(\X(n))$ is two, that 
$\sigma_{2n+2}(\X(n))$ is the entire envelopping projective space.

\bigskip
\bigskip
{\bf 3. Segre Veronese Varieties.} \par \medskip

Up to this point we have only considered the Segre varieties, i.e. the embeddings of $\P^{n_1}\times \cdots \times\P^{n_t}$ given by the very ample sheaves ${\cal O}(1, \ldots , 1)$.  We can also consider
the embeddings of these same varieties using the very ample sheaves ${\cal{O}}(d_1, \ldots , d_t)$, where $d_i > 0$.

These sheaves give a {\it Segre-Veronese} embedding (see \ref{BM} and \ref{CGG1}) into the projective space $\P^N$, where $N = (\prod_{i=1}^t N_i) -1$ and where $N_i = {n_i+d_i \choose n_i}$.  If we let $\underline{n} = (n_1, \ldots , n_t)$ and let $\underline{d} = (d_1, \ldots , d_t)$ then we will denote
this embedding by $\phi_{\underline{n}, \underline{d}}$ and its image by $\Xnd$.

If we denote by $\nu_{n_i,d_i}$ (or simply by $\nu_{d_i}$, when no doubt can occur), the Veronese embedding of $\P^{n_i}$ using forms of degree $d_i$, then $\phi_{\underline{n}, \underline{d}}$ is nothing
more than the composition:
$$
\psi \circ (\nu_{n_1,d_1}, \cdots , \nu_{n_t,d_t}) :\ \ \P^{n_1}\times\cdots\times\P^{n_t}\ \  \rightarrow \ \ \P^{N_1}\times\cdots \times\P^{N_t}\rightarrow \P^N
$$
where $\psi$ is just the usual Segre map and the other map is simply the product of the various Veronese
embeddings. 

To simplify the notation we just write:
$$
\P^{n_1}\times\cdots\times\P^{n_t}\ \buildrel{(d_1, \ldots , d_t)}\over{\longrightarrow} \P^N.
$$

In particular, we have that for ${\underline d}=(1,...,1)$, \ $\X _{{\underline n},(1,...,1)}$ is a Segre variety, while, for $t=1$, $\X_{\underline{n},\underline{d}}$ is the Veronese variety $\nu_{n,d}(\PP n)$.

Let $M_i$ be the $(n_i+1)\times {n_i+d_i-1 \choose n_i}$ catalecticant matrix whose $2\times 2$ minors define the ideal of the Veronese embedding $\nu_{d_i} (\PP {n_i})$ in $\PP {N_i}$.  Consider the matrix $M = M_1\otimes ... \otimes M_t$; \ since the $\nu_{d_i} (\PP {n_i})$'s are the rank 1 locus of $M_i$, $\Xnd$ is contained in the rank 1 locus of $M$.

Since, for generic matrices, the locus of the $(s+1)\times(s+1)$ minors is precisely the variety $\sigma_s(\Y)$, where $\Y$ is the locus of the $2\times 2$ minors, it follows that $\sigma _s(\Xnd)$ is contained in the zero locus of the $(s+1)\times (s+1)$ minors of $M$.  We get, in this way, equations for $\sigma _s(\Xnd)$. In this case, however, the matrix $M$ is not made up of independent coordinates (it has many equal entries, for example), and so we cannot ``a priori" know the heights of the ideals given by its minors just knowing their size. 

Nevertheless, whenever we expect $\sigma _s(\Xnd)$ to fill $\PP N$, and yet the $(s+1)\times (s+1)$ minors of $M$ give equations for $\sigma _s(\Xnd)$, we can definitely say that $\sigma _s(\Xnd)$ is defective.  To illustrate this consider the following two examples.
\par \bigskip

{\bf Example 3.1:}\ Consider $t=2$, $d_1=d_2=2$ and $n_1=n_2=n$, i.e.
$$
\PP n \times \PP n \quad \matrix {(2,2)\cr \longrightarrow}\quad \PP N, \ \ N = {n+2\choose 2}^2 - 1.
$$
Then for $n=1,2,3$ we have that $\Xnd$ is $s-defective$, for $s=n^2+2n$.

In fact, for these cases $M$ is an $(n+1)^2\times (n+1)^2$ matrix, and thus its determinant is zero on $\sigma _s(\Xnd)$, for $s=n^2+2n$.  Hence $\dim \sigma_{n^2 + 2n}(\X_{ (n,n), (2,2)})$ is $\leq N - 1$.  But, the expected dimension of $\sigma _s(\X_{(n,n),(2,2)})$ is $e=s(2n)+s-1$, and a straightforward computation shows that  $e \geq N$, for $n=1,2,3$.  

More precisely,

{\bf n=1:} \hskip 3cm{ $\PP 1 \times \PP 1 \quad \matrix {(2,2)\cr \longrightarrow}\quad \PP 8$ }

The image is the DelPezzo surface $D_8 \subset \P^8$ and will be discussed in detail in the next section.
 \bigskip

{\bf n=2:} \hskip 3cm{ $\PP 2 \times \PP 2 \quad \matrix {(2,2)\cr \longrightarrow}\quad \PP {35}$ }

Using, in a subtle manner, Horace's Method (see \ref{CGG1}) we obtain:  $\sigma_t(\X_{(2,2), (2,2)})$ has the expected dimension for $t \leq 6, \ t \geq 9$; $\sigma_8(\X_{(2,2), (2,2)})$ is a hypersurface whose equation is given above; the dimension of $\sigma_7(\X_{(2,2), (2,2)})$ is 32 rather than 34 (the expected dimension).  The $8\times 8$ minors of $M$ give us equations in the ideal of $\sigma_7(\X_{(2,2), (2,2)})$ but we do not know if they generate that ideal.
\bigskip

{\bf n=3:} \hskip 3cm{ $\PP 3 \times \PP 3 \quad \matrix {(2,2)\cr \longrightarrow}\quad \PP {99}$ }

Not only is  $\X_{(3,3), (2,2)}$ 14-defective (using the determinant of $M$ above) it is also 13-defective.  We conjecture that all the other secant varieties of $\X_{(3,3), (2,2)}$ have the expected dimension.

\bigskip

{\bf Example 3.2:}\ Consider $t=2$, $d_1=2k$, $d_2=2$, $n_1=1$ and $n_2=m$, i.e.  
$$
\PP 1 \times \PP m \quad \matrix {(2k,2)\cr \longrightarrow}\quad \PP N, N = (2k+1){m+2\choose 2} -1.
$$ 
Then $\forall m\geq 1, \quad k\geq 1$, we have that $\Xnd$ is $s-defective$, for $s=km+k+m$.  In fact, for these cases the $2k$-uple embedding of $\PP 1$ is defined by a $(k+1)\times (k+1)$ matrix $M_1$, while the 2-uple embedding of $\PP m$ is defined by $M_2$ of size $(m+1)\times (m+1)$; hence
$M=M_1\otimes M_2$ is an $(m+1)(k+1)\times (m+1)(k+1)$ matrix, and its determinant is zero on $\sigma _s(\Xnd)$, $s=mk+m+k$.

But(see \ref {CGG1, \ \S 3}), for such a variety the value $s_0=km+k+\lceil{m+1\over 2}\rceil\leq km+k+m$, is the one for which we expect that $\sigma _{s_0}(\Xnd)=\PP N$. Hence we have that $\sigma _s(\Xnd)$ is $s$-defective for all $s$ such that $s_0\leq s \leq km+k+m$.
 
\medskip

We'd like to point out that all the examples we found in \ref {CGG1, \S 3} can be viewed in this light, but this point of view also gives an equation for $\sigma _s(\Xnd)$).  Indeed, for these examples we are able to find a single (determinantal) equation to demonstrate that $\sigma _s(\Xnd)$ is defective. We
conjecture that this equation is {\bf the} defining equation for $\sigma _s(\Xnd)$
and, as a consequence, that $\sigma _{s+1}(\Xnd)=\PP N$.
 
\bigskip

\noindent {\bf 4. A particular case:  Del Pezzo Surfaces.}

Here we want to investigate the ideal of the classically studied Del Pezzo surfaces and of their secant varieties. Let $S_9$,$S_8$,...,$S_3$,$D_8$ be the (smooth) Del Pezzo surfaces of degree $d$ in $\PP d$, where $S_i\subset \PP i$ is obtained by blowing up $\PP 2$ at $9-i$ generic points, $i=3,...,9$ and then
embedding this into $\PP i$ via the linear system given by the strict transforms of the cubic curves passing through the points.   $D_8\subset \PP 8$, instead, is given by the embedding of a smooth
quadric $Q\subset \PP 3$ via the linear system given by ${\cal O}_Q(2)$, i.e. $D_8$ is the Segre-Veronese variety given by $\PP 1\times \PP 1$ embedded in $\P^8$ via ${\cal O}(2,2)$.

Even though these varieties have been extensively studied in classical Algebraic Geometry, their ideals (or at least some of them) seem not to be widely known.  For lack of a better reference we describe
them here. All of them, except $S_3$, have ideals which are generated by quadrics ($S_3$ is a surface in $\PP 3$ whose equation is the determinant of a $3\times 3$ matrix of linear forms).

The ideal of $S_4$ is generated by two quadrics and the ideal of $S_5$ by three quadrics. To see how these equations can be obtained from the generators of the ideal of the points in $\PP 2$ which
have been blown up see e.g. \ref {GiLo}.

The ideal of $S_9$, which is the 3-uple (Veronese) embedding of $\PP 2$ is well known; let $k[y_{000},...,y_{222}]$, $i,j,k \in \{0,1,2\}$, $i\leq j\leq k$ be the coordinate ring of $\PP 9$ and let the embedding $\PP 2 \rightarrow \PP 9$ be the morphism associated to the map $\kappa [y_{000},...,y_{222}] \rightarrow R = k[x_0,x_1,x_2]$ such that $y_{ijk}\rightarrow x_ix_jx_k$. Then $I_{S_9}$ is generated by the $2\times 2$ minors of the catalecticant matrix, $A_0$, which describes the multiplication $R_1 \times R_2 \rightarrow R_3$:
$$
A_0 = \pmatrix {y_{000}&y_{001}&y_{002}&y_{011}&y_{012}&y_{022}\cr
y_{001}&y_{011}&y_{012}&y_{111}&y_{112}&y_{122}
 \cr y_{002}&y_{012}&y_{022}&y_{112}&y_{122}&y_{222}}.
 $$

The ideals of $S_8$, $S_7$ and $S_6$ are known to be generated by quadrics (e.g. see \ref {Gi}); we will check that these quadrics can be obtained as $2\times 2$ minors of the matrix above just by erasing the last, then the fourth and then the first column of $A_0$.

In order to see why this is true we can view $S_8$ as given by the linear system of plane cubics containing (0:0:1), i.e. those cubics whose defining equations does not contain the monomial $x_2^3$; hence $S_8$ is the projection of $S_9$ onto the $\PP 8$ given by $y_{222}=0$ from the point $(0:...:0:1)\in S_9$. Then we have that $I_{S_8}= I_{S_9}\cap\kappa [y_{000},...,y_{122}]$.  Since we know that $I_{S_8}$ is generated in degree 2, its generators will be all the quadrics which are zero on $S_9$ and do not involve $y_{222}$. Those can be obtained by considering the $2\times 2$ minors of the matrix $A_1$ obtained by erasing from $A_0$ the column containing $y_{222}$, and the trick is done! Actually, all minors of $A$ involving the other two elements of that column are already given by the minors of $A_1$ and no linear combination of those involving $y_{222}$ gives new quadrics in $\kappa [y_{000},...,y_{122}]$.

In the same way we get the matrices $A_2$ and $A_3$ ($A_2$ by erasing the fourth column from $A_1$ and then $A_3$ by erasing the first column from $A_2$) whose $2\times 2$ minors give the ideals of $S_7$ and $S_6$.  All this can also be easily checked by \ref{CoCoA}.

Notice that those determinantal ideals, except for that of $S_6$ (see also \ref {GG}) are not generic, in the sense that they do not have the same height as the ideal of minors of a generic matrix of that size.

As for the ideal of $D_8$, working as in Example 3.1, we can see that it is generated by the $2\times 2$ minors of the $4\times 4$ matrix $B$ (this too can also be easily checked via \ref{CoCoA}):
$$
B = \pmatrix {y_{0000}&y_{0001}&y_{0100}&y_{0101}\cr
y_{0001}&y_{0011}& y_{0101}&y_{0111}\cr
y_{0100}&y_{0101}&y_{1100}&y_{1101}
 \cr y_{0101}&y_{0111}&y_{1101}&y_{1111}}
 $$

Here the embedding $\PP 1\times \PP 1 \rightarrow \PP 8$ is the one associated to the map
$$
k[y_{0000},...,y_{1111}] \rightarrow k[s_0,s_1,:t_0,t_1]; \quad y_{ijkl}\rightarrow
s_is_jt_kt_l,\quad \forall i,j,k,l \in \{0,1\},\ i\leq j;\  k\leq l.
$$
\medskip
 Now consider the varieties $\sigma_2(S_i)$, $i=3,...,9$. By Terracini's Lemma, we know that they have the expected dimension if the linear system of cubics passing through $9-i$ points and two double points (all generic) have the expected dimension (i.e. $max\{0,9-(9-i)-6\}$).  It is well-known (and easy to see) that this actually happens.  So, all the $S_i$'s are not 1-defective.

Also for $S_3$,$S_4$ and $S_5$ we have $\sigma_2(S_i)=\PP i$.  So, there is nothing to say about the ideal of $\sigma_2(S_i)$, $i=3,4,5$.

When  $i=6,7,8,9$ we want to show that the ideal of $\sigma _2(S_i)$ is generated by the $3\times 3$ minors of the matrices $A_{9-i}$ above.

First observe that, by a result of Kanev (see \ref {K}),  the ideal of $\sigma _2(S_9)$ is given by the $3\times 3 $ minors of $A_0$. Now consider the following remark:

\medskip \noindent {\bf Remark  4.1}  \  ``The secant variety of a projection is the projection of the
secant variety".   Let $X\subset \PP n$ be a non-degenerate reduced and irreducible variety, and $P\in \PP n$; let $\pi_P : \PP n - \{P\} \rightarrow H$ be the projection from $P$ to a generic hyperplane $H\cong \PP {n-1}$ and $X' = \overline{\pi_P(X - P)}$.  Then $\sigma_2(X')= \overline{\pi_P (\sigma_2(X)-P)}$.

\medskip
In fact, the inclusion $\sigma_2(X')\subset \overline{\pi_P (\sigma_2(X)-P)}$ is obvious.  As for the other inclusion, let $Q$ be a generic point of $\pi_P (\sigma_2(X)-P)$, and $Q'\in \sigma_2(X)-P$ a point in its preimage: there will be a secant line $L$ to $X$, not containing $P$ by genericity, which contains $Q'$, hence $\pi_P(L)$ will be a secant line of $\pi_P(X - P)$, and $Q \in \sigma_2(X')$.

\medskip Now, since $S_8,S_7$ and $S_6$ are obtained, each from the previous, by projection from one point at a time (starting from $S_9$), the same is true for their secant varieties (Remark 4.1).  Using 
 elimination, as we did for the ideals of $S_8,S_7, S_6$ themselves, we see that the ideals of $\sigma_2(S_i), i=6,7,8$ are given by the $3\times 3$ minors of the matrices $A_{9-i}$ (again, one can check this using \ref {CoCoA}).

Notice that in this case, i.e. for $\sigma_2(S_8),\sigma_2(S_7), \sigma_2(S_6)$ those determinantal ideals have generic height, hence are arithmetically Cohen-Macaulay with known resolution given by the Eagon-Northcott complex.

For $\sigma _2(D_8)$ we have to consider the $(2,2)$ divisors
through two generic 2-fat points in $\PP 1\times \PP 1$. It is
easy to check that this linear system has the expected dimension
(= 3), hence $\sigma _2(D_8)$ has the expected dimension (=5).

Is the ideal of $\sigma _2(D_8)$ generated by the $3\times 3$
minors of $B$?. By using \ref {CoCoA} one can check that this is
the case (notice that it does not have generic height).

From \ref {CS} we know the degree of $\sigma _2(D_8)$: $\deg \sigma _2(D_8)=10$.  Notice also that the degree of $\sigma _2(S_8)$ is 10.  A quick and easy check with {\bf CoCoA} shows that all the graded Betti numbers of these two varieties are also equal.  Hence up to this point $S_8$ and $D_8$ (and also their chordal varieties) cannot be distinguished by numerical invariants.

As for $\sigma _2(S_9)$, its degree is 15 (by \ref {CS} again) and this agrees with the fact that its ideal is generated  by the $3\times 3$ minors of $A_0$ (this was known by \ref K).

From Example 3.1 we know that $\sigma _3(D_8)$ does not have the expected dimension. It should fill $\PP 8$, but instead it is a hypersurface (see \ref{CGG1}).  Moreover, the equation of $\sigma _3(D_8)$ is given by $\det B$.

This is an interesting difference between $D_8$ and $S_8$: $\sigma _3(S_8)$ fills up $\PP 8 $ as expected, since there are no cubics in $\PP 2$ passing through three double points and a simple one.

In the same way we get that $\sigma _3(S_i)=\PP i$ for $i=6,7$, instead $\sigma _3(S_9)$ is a  hypersurface, as expected. Actually, $\sigma _3(S_9)$ is the hypersuface parameterizing Fermat cubics, so its equation (of degree 4) is defined by the Aronhold (or Clebsh) invariant of a cubic (e.g. see \ref {Ge} or \ref {DK}).

\bigskip \bigskip

\centerline {{\bf REFERENCES}}

\medskip\noindent[{\bf AOP}]: H. Abo, G. Ottaviani, C. Peterson. {\it Induction for secant varieites of Segre varieties,}\ arXiv:math.AG${\backslash}$ 0607191 v2 18 July 2006.

\medskip\noindent \ref{AR}:E.S. Allman, J.A. Rhodes. {\it Phylogenetic ideals and varieties for the general Markov model,}\ arXiv:math.AG${\backslash}$ 0410604 28 October 2004.

\medskip\noindent [{\bf BC}]: W. Bruns, A. Conca {\it Gr\"{o}bner Bases and Determinantal Ideals}, Commutative Algebra, Singularities and Computer Algebra (2003) Eds. J. Herzog, V. Vuletscu, Kluwer Academic Publishers, The Netherlands pp. 9-66.

\medskip\noindent [{\bf BCS}]: P. B\"{u}rgisser, M. Clausen, M.A. Shokrollahi. {\it Algebraic Complexity Theory.}, Vol. 315, Grund. der Math. Wiss., Springer, 1997.

\medskip\noindent [{\bf BM}]: S¸.Barcanescu, N. Manolache. Betti numbers of Segre-Veronese singularities. Rev. Roumaine Math. Pures Appl. {\bf 26} (1981), 549–565.

\medskip\noindent [{\bf CGG1}]: M.V.Catalisano, A.V.Geramita, A.Gimigliano. {\it Higher Secant Varieties of
Segre-Veronese varieties.} In: Varieties with unexpected properties. Siena,  Giugno 2004. BERLIN: W. de Gruyter. (2005), 81 - 107.

\medskip\noindent [{\bf CGG2}]: M.V.Catalisano, A.V.Geramita, A.Gimigliano. {\it Ranks of Tensors, secant varieties of Segre varieties and fat points.} Lin.Alg. and its Applications {\bf 355} (2002), 263-285.
(see also the errata of the publisher: {\bf 367} (2003) 347-348).

\medskip\noindent [{\bf CGG3}]: M.V.Catalisano, A.V.Geramita, A.Gimigliano. {\it Higher Secant varieties of the Segre varieties $\PP 1 \times ...\times \PP 1$}.  Jo. of Pure and Appl. Alg., (special volume in honour of W. Vasconcelos), vol. 201 (2005), 367-380.
 
\medskip\noindent [{\bf ChCi}]: L.Chiantini, C.Ciliberto. {\it Weakly defective varieties.} Trans. Am. Math. Soc. {\bf 354} (2001), 151-178.

\medskip\noindent [{\bf ChCo}]: L.Chiantini, M.Coppens. {\it Grassmannians for secant varieties.} Forum Math. {\bf 13} (2001), 615-628.

\medskip\noindent {\bf CoCoA}: A. Capani, G. Niesi, L. Robbiano, {\it CoCoA, a system for doing Computations in Commutative Algebra}  (Available via anonymous ftp from: cocoa.dima.unige.it).

\medskip\noindent \ref {CS}: D.Cox. J. Sidman. {\it Secant Varieties of toric varieties.} Preprint 2005. arXiv:math.AG/0502344
 
\medskip\noindent [{\bf DF}]: C. Dionisi, C.Fontanari. {\it Grassmann defectivity \`a la Terracini.} Le Matematiche {\bf 56}, (2001), 245-255.
 
\medskip\noindent \ref {DK}: I.Dolgachev, V. Kanev {\it Polar covariants of cubics and quartics.} Advances in Math., {\bf 98} (1993), 216-301.
 
\medskip\noindent [{\bf GSS}]: L.D. Garcia, M. Stillman, B. Sturmfels. {\it Algebraic Geometry of Bayseian Networks}, preprint 2003.

\medskip\noindent [{\bf GHKM}]: D.Geiger, D.Hackerman, H.King, C.Meek.  {\it Stratified Exponential Families: Graphical Models and Model Selection.} Annals of Statistics, {\bf 29} (2001), 505-527.

\medskip\noindent [{\bf Ge}]:  A.V.Geramita. {\it Inverse Systems of Fat Points: Waring's Problem, Secant
Varieties of Veronese Varieties and Parameter Spaces for Gorenstein Ideals.} The Curves Seminar at Queen's, Vol. X , Edited by A. V. Geramita, Queen's Papers in Pure and Applied Mathematics, {\bf 102} (1996), 3-104.

\medskip\noindent [{\bf GG}]: A.V.Geramita, A.Gimigliano. {\it Generators for the defining ideal of certain rational surfaces.} Duke Math. Journal, {\bf 62}, (1991), 61-83.

\medskip\noindent [{\bf Gi}]: A.Gimigliano. {\it On Veronesean Surfaces}. Indagationes Math., Ser.A, {\bf 92}  (1989).

\medskip\noindent[{\bf H}]:  J. Harris. {\it Algebraic Geometry: A First Course}. Springer-Verlag, 1992.

\medskip\noindent [{\bf GiLo}]: A.Gimigliano, A.Lorenzini. {\it On the ideal of some Veronesean Surfaces.} Canadian J.Math. {\bf 45}, (1993), 758-777.

\medskip\noindent [{\bf J}]:  J. Ja'Ja.\  {\it Optimal evaluation of pairs of bilinear forms.}, SIAM J. Comput. {\bf 8}, (1979), 443-462. 

\medskip\noindent [{\bf JPW}]:  T. Jozefiak, P. Pragacz and J. Weyman, {\it Resolutions of determinantal varieties and tensor complexes associated with symmetric and antisymmetric matrices}, Asterisque 87-88 (1981), 109--189. 

\medskip\noindent [{\bf K}]: V.Kanev. {\it Chordal varieties of Veronese varieties and catalecticant
 matrices.} J. Math. Sci.  {\bf 94}, (1999), 1114-1125.

\medskip\noindent [{\bf L}]: A. Lascoux { \it Syzygies des varietes d\'{e}terminatales (in French)}, Adv. in Math. {\bf 30} (1978), no. 3, 202-237.

\medskip\noindent [{\bf LM}]: Landsberg J.M., Manivel L. {\it On the ideals of secant varieties of Segre
 varieties.} Preprint, math.AG/0311388

\medskip\noindent [{\bf LW}]:  Landsberg J.M., Weyman J. {\it On the ideals and singularities of secant varieties of Segre varieties.}   Preprint, math.AG/0601452

\medskip\noindent [{\bf Pa}]: F. Palatini. {\it Sulle variet\`a algebriche per le quali sono di
dimensione minore} {\it dell' ordinario, senza riempire lo spazio ambiente, una o alcuna delle
variet\`a} {\it formate da spazi seganti.} Atti Accad. Torino Cl. Scienze Mat. Fis. Nat. {\bf 44}
(1909), 362-375.

\medskip\noindent [{\bf P}]: Parolin, A., {\it Variet\`a delle secanti di variet\`a di Segre e Veronese e Applicazioni.}  ({\it Tesi di Dottorato}), 2004, Universit\`a di Bologna, Italy

\medskip\noindent [{\bf Te}]: A. Terracini. {\it Sulle} $V_k$ {\it per cui la variet\`a degli} $S_h$ $(h+1)${\it -seganti ha dimensione minore dell'ordinario.} Rend. Circ. Mat. Palermo {\bf 31} (1911), 392-396.

\medskip\noindent [{\bf Te2}]: A. Terracini. {\it Sulla rappresentazione delle coppie di forme ternarie mediante somme di potenze di forme lineari.} Ann. Mat. Pur ed appl. {\bf XXIV}, {\bf III} (1915),
91-100.

\medskip\noindent[{\bf W}]:  Weyman, J.{\it Cohomology of Vector Bundles and Syzygies} ISBN: 0511059701,  Cambridge University Press, 2003

\medskip\noindent [{\bf Z}]: F.L.Zak. {\it Tangents and Secants of Algebraic Varieties.} Translations of Math. Monographs, vol. {\bf 127} AMS. Providence (1993).

\bigskip

{\it M.V.Catalisano, DIPTEM, Univ. di Genova, Italy.}

{\it e-mail: catalisano@diptem.unige.it}

\medskip
{\it A.V.Geramita, Dept. Math. and Stats. Queens' Univ. Kingston, Canada}
{\it and Dip. di Matematica, Univ. di Genova. Italy.}

{\it e-mail: geramita@dima.unige.it ; tony@mast.queensu.ca}

\medskip
{\it A.Gimigliano, Dip. di Matematica and C.I.R.A.M., Univ. di Bologna, Italy.}

{\it e-mail: gimiglia@dm.unibo.it}

\end